\documentclass{article}

\usepackage[latin1]{inputenc}
\usepackage{latexsym}
\usepackage{amssymb, amsmath}
\usepackage{amsfonts}
\usepackage{epsfig}
\usepackage{color}
\usepackage{graphicx,graphics}

\newtheorem{prop}{Proposition}[section]

\newtheorem{lemma}{Lemma}[section]
\newtheorem{cor}{Corollary}[section]

\newcommand{\cvd}{\hfill $\blacksquare$\bigskip}

\author{Filippo Disanto\thanks{Dipartimento di Scienze Matematiche ed
Informatiche, Pian dei Mantellini, 44, 53100, Siena, Italy \; {\tt
disafili@yahoo.it\quad rinaldi@unisi.it}}\and Luca
Ferrari\thanks{Dipartimento di Sistemi e Informatica, viale
Morgagni 65, 50134 Firenze, Italy {\tt ferrari@dsi.unifi.it\quad
pinzani@dsi.unifi.it}}\and Renzo Pinzani$^\dag$ \and Simone
Rinaldi$^*$}

\title{Catalan lattices on series parallel interval orders}

\begin{document}

\date{}

\maketitle

\begin{abstract}
Using the notion of series parallel interval order, we propose a
unified setting to describe Dyck lattices and Tamari lattices (two
well known lattice structures on Catalan objects) in terms of
basic notions of the theory of posets. As a consequence of our
approach, we find an extremely simple proof of the fact that the
Dyck order is a refinement of the Tamari one. Moreover, we provide
a description of both the weak and the strong Bruhat order on
312-avoiding permutations, by recovering the proof of the fact
that they are isomorphic to the Tamari and the Dyck order,
respectively; our proof, which simplifies the existing ones,
relies on our results on series parallel interval orders.
\end{abstract}

\section{Introduction}

\begin{figure}[htb]
\begin{center}
\centerline{\hbox{\includegraphics[width=2.8in]{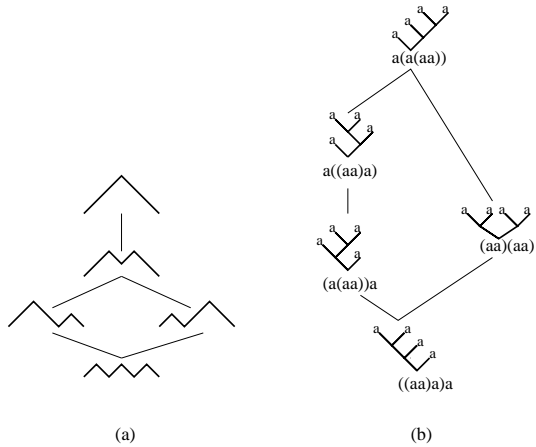}}}
\caption{(a) The Dyck lattice with five elements; (b) The Tamari
lattice with five elements.} \label{latt}
\end{center}
\end{figure}

The set of all Dyck paths of fixed length  can be ordered by
setting $P\leq Q$ when $P$ lies weakly below $Q$ (in the usual
two-dimensional drawings of Dyck paths). This partial order is in
fact a lattice, called \emph{Dyck lattice}. We point out that this
structure is not new and, in some sources, it is called
\emph{Stanley lattice} \cite{BB,K4}. The Hasse diagram of the Dyck
lattice on the set of Dyck paths of length 6 is represented in
figure \ref{latt}(a).

It is well known that Dyck paths of length $2n$ are counted by the
$n$-th \emph{Catalan} number $C_n =\frac{1}{n+1}$${2n}\choose{n}$.
The different incarnations of the Catalan family give rise to
several further lattices beside Dyck's. Among them, the Tamari
lattice is indeed one of the more widely known, and appears
naturally in the study of binary trees and of the Stasheff
polytope \cite{HT}. See figure \ref{latt}(b) for the Hasse diagram
of the Tamari lattice with five elements.



Using suitable bijections between Dyck paths, binary trees and
planar trees, the two mentioned Catalan lattices can be defined on
the set of plane trees of size $n$ in such a way that the Dyck
lattice with $n$ elements is an extension of the Tamari lattice
with $n$ elements (see \cite{K4}).

In this paper, we will consider yet another occurrence of Catalan
structures, the so called \emph{series parallel interval orders}.
Our aim is to define the two above Catalan lattices on the set of
series parallel interval orders with the aid of a special kind of
linear extension. We propose this unified interpretation since, in
our opinion, it allows to better understand the connections
between Dyck and Tamari lattices. To obtain our characterization
of the Dyck and Tamari lattices we will make use of some basic
notions of the theory of posets. The known relationship between
the two lattices will be quite simple to achieve in our setting.

Thanks to this approach we will also be able to provide a link
between the Dyck (respectively, Tamari) lattice and the strong
(respectively, weak) Bruhat order, when the latter is considered
on the class of permutations avoiding the pattern 312.

As already recalled, the main combinatorial objects in this
approach are series parallel interval orders, which are the
intersection of two important classes of partially ordered sets,
namely series parallel orders \cite{S1} and interval orders
\cite{B-MCDK,F}. Our approach will express important features of
series parallel interval orders and so their use in this unified
version of the two Catalan lattices seems to be relevant in its
own.

\section{Series-parallel interval orders}\label{spio}

\begin{figure}[htb]
\begin{center}
\centerline{\hbox{\includegraphics[width=1.8in]{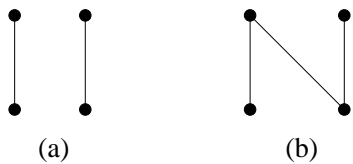}}}
\caption{(a) The poset $2+2$; (b) The fence of order four.}
\label{poset}
\end{center}
\end{figure}

In this section we will focus on those posets having no induced
subposet isomorphic either to the poset $2+2$ or to the fence of
order four, shown in figure \ref{poset}. These partial orders are
called \emph{series parallel interval orders}. We will denote by
$\cal{O}$ the class of such posets, also writing ${\cal{O}}(n)$
for those having precisely $n$ elements. We also warn the reader
that, due to technical reasons, in what follows we will rather
deal with the strict order relation associated with a series
parallel interval order. Nevertheless, with an abuse of notation,
we will always use the expression $R\in \mathcal{O}$ to mean that
$R$ is the strict order relation associated with a series parallel
interval order.

\begin{figure}[htb]
\begin{center}
\centerline{\hbox{\includegraphics[width=3.8in]{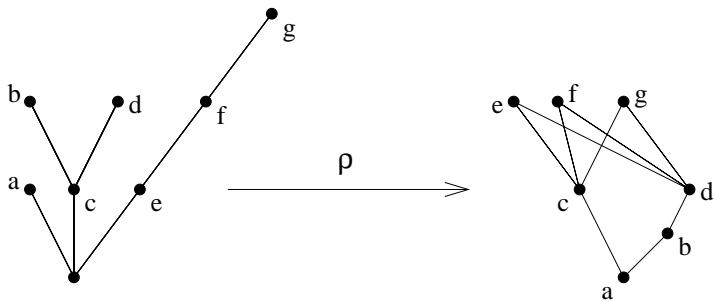}}}
\caption{The bijection $\rho$: a planar tree and the associated
poset in $\cal{O}$.} \label{tree}
\end{center}
\end{figure}

This kind of posets has been recently considered in \cite{DFPR},
where the authors show that they are enumerated by Catalan numbers
according to the number of their elements; some bijections with
other structures enumerated by Catalan numbers are also
established. For our purposes, we need to recall here a bijection
$\rho$ (stated in \cite{DFPR}) between planar trees with $n+1$
nodes and $\mathcal{O}(n)$. Given any planar tree $T$, we define
the binary relation $R=\rho(T)$ on the set of its nodes other than
the root, by setting $xRy$ whenever $x$ and $y$ cannot be joined
by a directed path in $T$ (in the directed graph canonically
determined by $T$) and $x$ lies on the left of $y$ in $T$. The
resulting poset is indeed in $\mathcal{O}$. In figure \ref{tree}
we can see an instance of the bijection $\rho$. In what follows,
we will always represent rooted trees with the root at the bottom.

\subsection{Preorder linear extensions of series parallel interval orders}

In this section we will define a particular type of linear
extension for the posets in $\cal{O}$ corresponding to the
preorder traversal in the associated planar tree; for this reason
we will call it the \emph{preorder linear extension}.

%
%

\begin{figure}[htb]
\begin{center}
\centerline{\hbox{\includegraphics[width=2.8in]{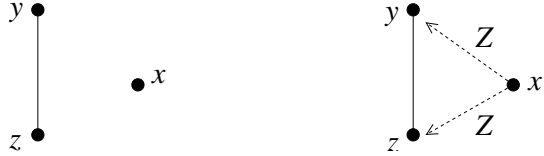}}}
\caption{The relation $Z$ on an instance of the poset $2+1$.}
\label{2+1}
\end{center}
\end{figure}

In order to define the preorder linear extension of $R\in
\cal{O}$, we need to define an auxiliary binary relation $Z(R)=Z$
on the support of $R$. Given a binary relation $B$, we set
$\overline{B} = B \cup B^{-1}$ and we use the notation $B^c$ to
indicate the complement of $B$. Now define $Z = ((\overline{R})^c
\circ \overline{R}) \setminus \overline{R}$. Recall that, for any
two binary relations $X$ and $Y$ defined on the same set, the
\emph{composition} $X\circ Y$ is defined by setting $x(X\circ Y)y$
when there exists an element $z$ such that $xXz$ and $zYy$ (see,
for instance, \cite{KRY}). Thus, we can rephrase the above
definition by saying that $xZy$ whenever $x\! \! \not \! \!
\overline{R}y$ and there exists $z$ such that $z\overline{R}y$ and
$z\! \! \not \! \! \overline{R}x$. Given $R\in \cal{O}$, if
$x,y,z$ are such that $zRy$ and $x$ is incomparable with both $z$
and $y$, then $Z$ can be described as illustrated in figure
\ref{2+1}.

\begin{figure}[htb]
\begin{center}
\centerline{\hbox{\includegraphics[width=2.8in]{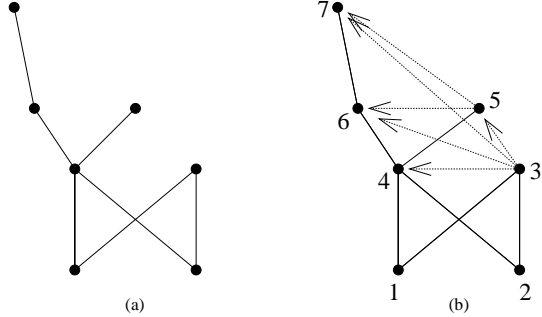}}}
\caption{(a) An element $R \in \cal{O}$; (b) The relation $Z$
associated with $R$, where pairs of $Z$ are joined with an arrow.
Vertices are then labelled according to a preorder linear
extension.} \label{extension}
\end{center}
\end{figure}

We say that a linear extension $\lambda$ of $R$ is a
\emph{preorder linear extension} of $R$ when $xZy$ implies
$x\lambda y$. Figure \ref{extension} depicts a poset $R \in
\cal{O}$, the relation $Z(R)$ and an associated preorder linear
extension.

The next proposition shows that, for any $R \in \cal{O}$, there
exists at most one preorder linear extension of $R$ up to (order)
automorphisms. The proposition needs a preliminary lemma. In what
follows, we will say that two elements of a poset are \emph{order
equivalent} when there exists an order automorphism mapping one of
them into the other.

\begin{lemma}
Let $R\in \cal{O}$ and suppose that $\lambda_1$ and $\lambda_2$
are two preorder linear extensions of $R$. For any $x,y$ in the
support of $R$, if $\lambda_1(x)>\lambda_1(y)$ and
$\lambda_2(x)<\lambda_2(y)$, then $x$ and $y$ must be order
equivalent in $R$.
\end{lemma}

\emph{Proof.}\quad If $\lambda_1(x)>\lambda_1(y)$ then $(x,y)
\notin R\cup Z$ and by $\lambda_2(x)<\lambda_2(y)$ follows that
$(y,x) \notin R\cup Z$. From the definition of $Z$, this implies
that, for every $z$, $z\! \not \! \! \overline{R}y$ or
$z\overline{R}x$, and also that, for every $z$, $z\overline{R}y$
or $z\! \not \! \! \overline{R}x$. Equivalently, for every $z$, we
get that either $z\! \not \! \! \overline{R}y$ and $z\! \not \! \!
\overline{R}x$, or $z \overline{R}y$ and $z \overline{R}x$. Since
$x\! \not \! \! \overline{R}y$, it is now easy to show that, for
every $a\neq x,y$, $aRx$ if and only if $aRy$ and $xRa$ if and
only if $yRa$, whence the thesis.\cvd

\begin{prop} \label{unique}
Let $R\in \cal{O}$ and suppose that $\lambda_1$ and $\lambda_2$
are two preorder linear extensions of $R$. If $\lambda_1(x) =
\lambda_2(y)$, then the two elements $x$ and $y$ must be order
equivalent in $R$.
\end{prop}

\emph{Proof.}\quad If $\lambda_1(x)>\lambda_1(y)$ and
$\lambda_2(x)<\lambda_2(y)$, then the thesis follows from the
above lemma. Otherwise, without loss of generality, suppose that
$\lambda_1(x)>\lambda_1(y)$ and $\lambda_2(x)>\lambda_2(y)$. We
claim that there exists an element $\tilde{z}$ such that $x$ is
equivalent to $\tilde{z}$ which in turn is equivalent to $y$, and
this will be enough to conclude.

Indeed, denoting with $X$ the support of $R$, we have that $|\{z
\in X \setminus \{x,y\} : \lambda_2(z)> \lambda_2(y)=
\lambda_1(x)\}|+1=|\{z \in X \setminus \{x,y\} : \lambda_1(z)>
\lambda_2(y)= \lambda_1(x)\}|$. Then there must exist an element
$\tilde{z} \in X \setminus \{x,y\}$ such that
$\lambda_1(\tilde{z})>\lambda_1(x)=\lambda_2(y)$ and
$\lambda_2(\tilde{z})<\lambda_1(x)=\lambda_2(y)$. Since
$\lambda_1(\tilde{z})>\lambda_1(x)$ and
$\lambda_2(\tilde{z})<\lambda_2(y)<\lambda_2(x)$, then $\tilde{z}$
must be order equivalent to $x$ in $R$ (once again thanks to the
above lemma). Analogously, since
$\lambda_1(\tilde{z})>\lambda_1(x)>\lambda_1(y)$ and
$\lambda_2(\tilde{z})<\lambda_2(y)$, $\tilde{z}$ and $y$ are order
equivalent.\cvd

\bigskip

Thanks to the above proposition we can assert that, for any $R\in
\mathcal{O}$, there is at most one preorder linear extension of
$R$. The next proposition shows that indeed a (the) preorder
linear extension exists, and also suggests how to find it.

\begin{prop} Let $T$ be a planar tree and $\rho$ be the
bijection described in section \ref{spio}. Suppose that the nodes
of $T$ are labelled according to the preorder traversal. Then the
induced labelling on $\rho (T)$ determines a preorder linear
extension of $\rho (T)$.
\end{prop}

\emph{Proof.}\quad With a slight abuse of notation, in this proof
we will denote the elements of $T$ and $\rho (T)$ using their
labels in the appropriate linear extensions. Denote with $<_t$ the
total order determined by the preorder traversal on $T$. We have
to prove that $<_t$ is mapped by $\rho$ to a preorder linear
extension of $\rho (T)$. According to the definition of preorder
linear extension, what we have to show is that, given two nodes
$x$ and $y$ of $T$ such that $x<_t y$, then the pair $(x,y)$
satisfies the definition of preorder linear extension. More
precisely, we must prove that, if $xZy$ in $\rho (T)$, then
necessarily $x<_t y$. Indeed, $xZy$ immediately implies that $x\!
\! \not \! \! \overline{R}y$. This means that, in $T$, either $x$
is a descendant of $y$ or $y$ is a descendant of $x$. Suppose that
the former case holds. The fact that $xZy$ also implies that there
exists an element $z$ in $\rho (T)$ such that $z\overline{R}y$ and
$z\! \! \not \! \! \overline{R}x$. In particular, this would mean
that, in $T$, $x$ should be a descendant of $z$ and, at the same
time, neither $z$ could be a descendant of $y$ nor $y$ could be a
descendant of $z$, which is plainly impossible (since $T$ is a
tree). Therefore we must have that $y$ is a descendant of $x$ in
$T$, and so $x<_t y$.\cvd

\bigskip

We remark that the definition of preorder linear extension has a
meaning only for elements in $\cal{O}$: just observe that it is
not possible to construct a linear extension with the required
properties neither of 2+2 nor of the fence of order 4.

%
%

\begin{figure}[htb]
\begin{center}
\centerline{\hbox{\includegraphics[width=3.0in]{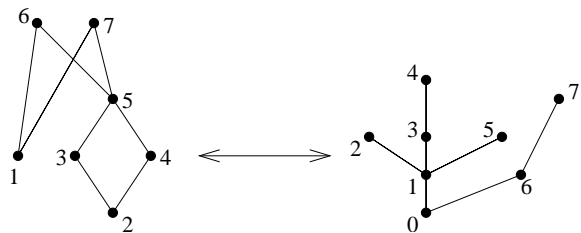}}}
\caption{A series parallel interval order with its preorder linear
extension and the associated planar tree with the labelling of its
node according to the preorder visit.} \label{equivalence1}
\end{center}
\end{figure}

In figure \ref{equivalence1} an example of the correspondence
between preorder traversal and preorder linear extension is shown.

\section{Catalan lattices on series parallel interval orders}

In this section we will define the Dyck lattice and the Tamari
lattice on series parallel interval orders whose support is
equipped with a preorder linear extension.

In the sequel we will refer to a node of a planar tree using its
label in the preorder traversal. Similarly, we will tacitly assume
that posets in $\cal{O}$ are equipped with their preorder linear
extension, and we will refer to their elements using the
corresponding labels. Moreover, if $x$ and $y$ are labels, we will
write $x<y$ referring to the usual order on natural numbers.

Given a planar tree $T$, $\rho(T)$ will denote the poset obtained
through the bijection $\rho$ defined above. In particular, we will
often refer to $\rho(T)$ as a binary relation. Given $x\in X$, for
a binary relation $B$ on $X$ we consider the set $B(x)=\{ y\in X\,
:xBy\}$. If $B$ is a partial order, then the set $B(x)$ is the
\emph{principal (order) filter} generated by the element $x$,
whereas $B^{-1}(x)=\{ y\in X : yBx \}$ is the \emph{principal
(order) ideal} generated by $x$. In the sequel, we will always use
the terms ``filter" and ``ideal" in place of ``order filter" and
``order ideal". Moreover, even if this is non standard in poset
theory, a principal filter (respectively, ideal) will be
considered without its minimum (respectively, maximum).

\subsection{The Dyck lattice}

According to \cite{K4}, we start by recalling the definition of
the Dyck lattice for the set of planar rooted trees with a fixed
number of nodes. If $T$ is a planar tree and $k$ is a node of $T$,
then define $h_T(k)$ as the set of ancestors of $k$ in the tree
$T$. Given two planar trees $T_1$ and $T_2$ having $n$ nodes,
$T_1$ is less than or equal to $T_2$ in the Dyck order, written as
$T_1 \leq_{D} T_2$, whenever, for every node $k$,
$|h_{T_{1}}(k)|\leq |h_{T_{2}}(k)|$.

The above definition allows us to give a characterization of the
Dyck order in terms of series parallel interval orders.

\begin{prop}
Let $T_1$, $T_2$ be two planar trees having $n$ nodes and let
$\rho(T_1) = R_1$ and $\rho(T_2) = R_2$. Then the following
conditions are equivalent:

\begin{itemize}
\item[a)] $T_1 \leq_{D} T_2$ ;

\item[b)] for every $k$, $|R_1 ^{-1}(k)|\geq |R_2 ^{-1}(k)|$.
\end{itemize}
\end{prop}

\emph{Proof.}\quad By definition $T_1 \leq_{D} T_2$ if and only
if, for all $k$, $|h_{T_{1}}(k)|\leq |h_{T_{2}}(k)|$, which is
equivalent to:
$$
|\{ x \in T_1:x\notin h_{T_{1}}(k) \}| \geq |\{ x \in
T_2:x\notin h_{T_{2}}(k)\} |.
$$

Consider now the series parallel interval orders $R_1$ and $R_2$
associated with $T_1$ and $T_2$ respectively. The previous
condition may be expressed by saying that, for all $k$:
\begin{equation}\label{temp}
|\overline{R_1}(k) \cup \{ x \in T_1:x>k \}| \geq
|\overline{R_2}(k) \cup \{ x \in T_2:x>k \}|.
\end{equation}

To show that (\ref{temp}) is equivalent to $b)$ observe that, for
a generic element $k$, the inequality
$$
|\overline{R_1}(k) \cup \{k+1,...,n\}| \geq |\overline{R_2}(k)
\cup \{k+1,...,n\}|
$$
holds if and only if
$$
|{R_1}^{-1}(k) \cup \{k+1,...,n\}| \geq |{R_2}^{-1}(k) \cup
\{k+1,...,n\}|,
$$
since, for $i=1,2$, $R_i (k)\subseteq \{ k+1,\ldots n\}$. Thus,
being also $R_i^{-1}(k)\cap \{ k+1\ldots, n\}=\emptyset$, we
immediately get
$$
|{R_1}^{-1}(k)| + |\{k+1,...,n\}| \geq |{R_2}^{-1}(k)| +
|\{k+1,...,n\}|,
$$
which is precisely $b)$.\cvd

\bigskip

\begin{figure}[htb]
\begin{center}
\centerline{\hbox{\includegraphics[width=2.8in]{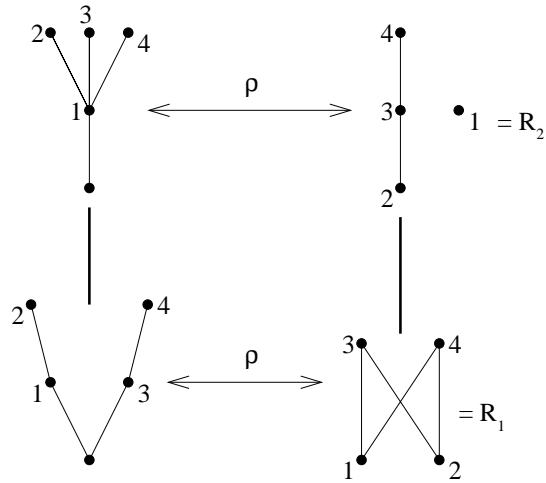}}}
\caption{Two comparable elements in the Dyck order.} \label{stan}
\end{center}
\end{figure}

In figure \ref{stan} an application of this proposition is shown.
The figure depicts two comparable elements in the Dyck order, and
the trees on the left correspond to the posets on the right
through $\rho$.

As a corollary, we find that the Dyck order on series parallel
interval orders can be defined in terms of the cardinalities of
the principal ideals of these posets.

\begin{cor}\label{ref1}
Given $R_1,R_2 \in \mathcal{O}(n)$, $R_1 \leq_{D} R_2$ if and only
if, for every $k\leq n$, $|R_1 ^{-1}(k)|\geq |R_2 ^{-1}(k)|$.
\end{cor}

See also figure \ref{fin} (a) for an example.

\subsection{The Tamari lattice}

According to \cite{K4}, we start by recalling the definition of
the Tamari lattice for the set of planar rooted trees with a fixed
number of nodes. If $T$ is a planar tree and $k$ is a node of $T$,
then define $u_T(k)$ as the set of descendants of $k$ in the tree
$T$. Given two planar trees $T_1$ and $T_2$ having $n$ nodes,
$T_1$ is less than or equal to $T_2$ in the Tamari order, written
as $T_1 \leq_{T} T_2$, whenever, for every node $k$,
$|u_{T_{1}}(k)|\leq |u_{T_{2}}(k)|$.

We provide two equivalent conditions to define the Tamari order on
the set of planar rooted trees.

\begin{lemma}
If $T_1 ,T_2$ are two planar trees with $n$ nodes, then $T_1
\leq_{T}T_2$ if and only if, for every node $k$, $u_{T_{1}}(k)
\subseteq u_{T_{2}}(k)$.
\end{lemma}

\emph{Proof.}\quad Given a node $k$ in the planar tree $T$, if
$|u(k)| = j$, then obviously $u(k)=\{
k+1,k+2,...,k+j-1,k+j\}$.\cvd

\bigskip

Given a node $k$ in the planar tree $T$, consider the set $d(k) =
u(k) \cup h(k)$, i.e. the set of ancestors and descendants of $k$.
The following lemma holds.

\begin{lemma}
If $T_1 ,T_2$ are two planar trees with $n$ nodes, then
$u_{T_{1}}(k)\subseteq u_{T_{2}}(k)$ (for every node $k$) if and
only if $d_{T_{1}}(k)\subseteq d_{T_{2}}(k)$ (for every node $k$).
\end{lemma}

\emph{Proof.}\quad If $u_{T_{1}}(k)\subseteq u_{T_{2}}(k)$ for all
$k$, then $h_{T_{1}}(k)\subseteq h_{T_{2}}(k)$ for all $k$ as
well. Indeed, if $x\in h_{T_{1}}(k)$, then $k \in u_{T_{1}}(x)$,
hence $k\in u_{T_{2}}(x)$ and so $x\in h_{T_{2}}(k)$. Therefore
$d_{T_{1}}(k)=u_{T_{1}}(k)\cup h_{T_{1}}(k)\subseteq
u_{T_{2}}(k)\cup h_{T_{2}}(k)=d_{T_{2}}(k)$, for every node $k$.

Now suppose that, for every node $k$, $d_{T_{1}}(k)\subseteq
d_{T_{2}}(k)$. If $x\in u_{T_{1}}(k)$ then $x\in d_{T_{1}}(k)$ and
then $x\in d_{T_{2}}(k)$ with $x>k$. So $x\in u_{T_{2}}(k)\cup
h_{T_{2}}(k)$, that is $x\in u_{T_{2}}(k)$, since $x>k$.\cvd

\bigskip

The above lemma allows us to give the following characterization
of the Tamari order.

\begin{prop}
Let $T_1$, $T_2$ be two planar trees having $n$ nodes and let
$\rho(T_1) = R_1$ and $\rho(T_2) = R_2$. Then the following
conditions are equivalent :

\begin{itemize}
\item[a)] $T_1 \leq_{T} T_2$;

\item[b)] for every $k$, $R_1 (k)\supseteq R_2 (k)$.

\end{itemize}

\end{prop}
\emph{Proof.}\quad From the previous lemma we have that $T_1
\leq_{T} T_2$ if and only if, for every $k$, $d_{T_{1}}(k)
\subseteq d_{T_{2}}(k)$. Now $d_{T_{1}}(k) \subseteq d_{T_{2}}(k)$
if and only if, for every $y$, $k R_2 y$ implies $k R_1 y$, that
is, for every $k$, $R_1(k) \supseteq R_2(k)$. Indeed, suppose
that, for every $k$, $R_1(k) \supseteq R_2(k)$. If $y \notin
d_{T_{2}}(k)$ then $k \overline{R_2} y$, whence $k \overline{R_1}
y$, and so $y \notin d_{T_{1}}(k)$. Viceversa, suppose that, for
every $k$, $d_{T_{1}}(k)\subseteq d_{T_{2}}(k)$. If
$y\overline{R_2}k$, then $y \overline{R_1}k$. Now, if $kR_2 y$,
then obviously $k<y$, which cannot hold together with $yR_1 k$ and
so $kR_1 y$. Thus we can conclude that $yR_2 k$ implies $yR_1 k$,
as desired.\cvd

\bigskip

\begin{figure}[htb]
\begin{center}
\centerline{\hbox{\includegraphics[width=2.8in]{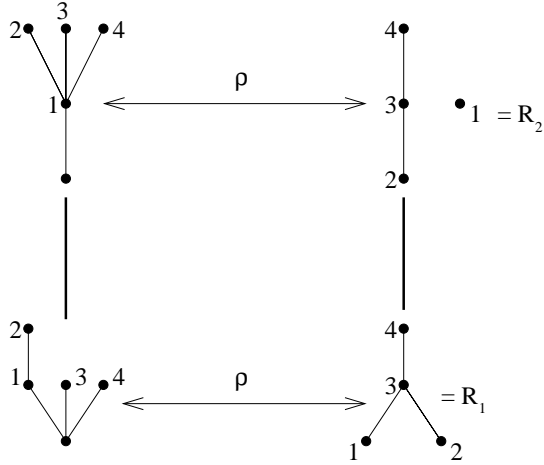}}}
\caption{Two comparable elements in the Tamari order.} \label{tam}
\end{center}
\end{figure}

Figure \ref{tam} shows an application of this proposition. The
figure depicts two comparable elements in the Tamari order, and
the trees on the left correspond to the posets on the right
through $\rho$.

\bigskip

The next lemma shows that the Tamari order may be defined by means
of the cardinalities of the principal filters of the posets in
$\mathcal{O}$.

\begin{lemma}
Let $T_1$, $T_2$ be two planar trees having $n$ nodes and
$\rho(T_1) = R_1$, $\rho(T_2) = R_2$ be the two associated posets
in $\mathcal{O}(n)$. Then, for every node $k$, $|R_1(k)| \geq
|R_2(k)|$ if and only if $R_1(k) \supseteq R_2(k)$.
\end{lemma}

\emph{Proof.}\quad Just observe that, if $T$ is a planar tree
having $n$ nodes, then, for any of its nodes $k$, $R(k)$ is a
final segment of $\{ 1,\ldots ,n\}$.\cvd

\bigskip

The above lemma allows us to give another characterization of the
Tamari order.

\begin{prop}
Let $T_1$, $T_2$ be two planar trees having $n$ nodes and let
$\rho(T_1) = R_1$ and $\rho(T_2) = R_2$. Then $T_1 \leq_{T} T_2$
if and only if, for every node $k$, $|R_1 (k)|\geq |R_2 (k)|$.
\end{prop}

\bigskip

As a corollary of the above propositions, we can immediately
obtain a definition of the Tamari order on series parallel
interval orders in terms of the principal filters of these posets.

\begin{cor}\label{ref2}
Given $R_1,R_2 \in \mathcal{O}(n)$, $R_1 \leq_T R_2$ if and only
if, for every $k\leq n$, $|R_1 (k)|\geq |R_2 (k)|$. Equivalently,
$R_1 \leq_T R_2$ if and only if, for every $k\leq n$, $R_1
(k)\supseteq R_2 (k)$.
\end{cor}

\begin{figure}[htb]
\begin{center}
\centerline{\hbox{\includegraphics[width=3.8in]{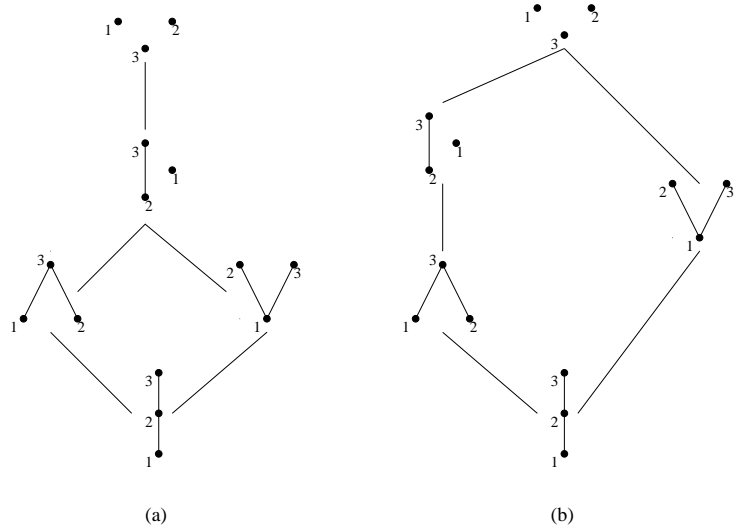}}}
\caption{The two Catalan lattices defined on the class of series
parallel interval orders of size three: (a)The Dyck lattice (b)The
Tamari lattice.} \label{fin}
\end{center}
\end{figure}

In figure \ref{fin}(b) the Tamari lattice on the five elements
belonging to ${\cal{O}}(3)$ is depicted .

\bigskip

\emph{Remark.} We know from \cite{K4} that the Dyck and the Tamari
orders are related by the following refinement property: given two
Catalan structures $T_1,T_2$ of the same size, if $T_1 \leq_T T_2$
then $T_1 \leq_D T_2$. This fact is an obvious consequence of
Corollaries \ref{ref1} and \ref{ref2}. Our approach seems to be
particulary interesting since it is now possible to prove such a
refinement property in a very neat way. Indeed, if $R_1
(k)\supseteq R_2 (k)$ holds for any $k$, then we also have that,
for all $k$, $R_1 ^{-1}(k)\supseteq R_2 ^{-1}(k)$, and so, for all
$k$, $|R_1 ^{-1}(k)|\geq |R_2 ^{-1}(k)|$.

\section{Series parallel interval orders and pattern avoiding permutations}

The \emph{strong Bruhat} order ($\leq_{B}$) and the \emph{weak
Bruhat} order ($\leq_{b}$) are two well known partial orders
defined on the set of permutations having fixed length \cite{S1}.
We briefly recall here their definitions.

Given a permutation $\pi = a_1a_2...a_n$, a \emph{reduction} of
$\pi$ is a permutation obtained from $\pi$ by interchanging some
$a_i$ with some $a_j$, provided that $i<j$ and $a_i > a_j$. We say
that $\pi_1 <_{B} \pi_2$ whenever $\pi_1$ is obtained from $\pi_2$
through a sequence of reductions. Define a \emph{simple reduction}
of $\pi = a_1a_2...a_n$ as a permutation obtained from $\pi$ by
interchanging some $a_i$ with some $a_{i+1}$, provided that $a_i >
a_{i+1}$. We say that $\pi_1 <_{b} \pi_2$ whenever $\pi_1$ is
obtained from $\pi_2$ through a sequence of simple reductions.

In this section we will consider another well known Catalan
structure, namely the class of permutations avoiding the pattern
$312$, and we will prove, using a characterization given in
\cite{BBFP}, that the strong Bruhat order, when restricted to such
a class of pattern avoiding permutations, is isomorphic to the
Dyck order. Moreover we will show that an analogous isomorphism
also exists between the Tamari lattice and the weak Bruhat order
on the same class of permutations. We remark that these two
results have been already obtained independently in \cite{BBFP}
(for the Dyck case) and in \cite{BW} (for the Tamari case; see
also \cite{D}). Here our main aim is to find a common language for
these two results.

We start by describing a bijection between series parallel
interval orders on $n$ elements and permutations of length $n$
avoiding the pattern $312$, denoted by $Av_n(312)$. Our approach
can be compared with the one used in \cite{B-MCDK} to enumerate
posets avoiding $2+2$.

First of all recall that the set of principal filters of a poset
avoiding 2+2 is linearly ordered by inclusion. The interested
reader can find a proof of this fact in \cite{E-Z}, where it is
also proved that this condition completely characterizes such a
class of posets.

Given a poset $R \in {\cal{O}}(n)$, consider the labelling of its
elements determined by its preorder linear extension and denote
its principal filters by $R(1),R(2),$ $\ldots ,$ $R(n)$. Define a
permutation $\pi = \pi(R)$ of length $n$ as follows: $k$ precedes
$j$ in $\pi$ precisely when either $R(k)\supset R(j)$ or
$R(k)=R(j)$ and $k>j$. It is easy to show that, for each $R \in
\cal{O}$, $\pi(R)$ does not contain the pattern 312, and the
function $R \mapsto \pi(R)$ is a bijection between ${\cal{O}}(n)$
and $Av_n(312)$.

\bigskip

\emph{Remark.} Observe that our bijection cannot be described in
terms of principal ideals (instead of principal filters), due to
our choice of taking the preorder linear extension of a poset.

\bigskip

\begin{figure}[htb]
\begin{center}
\centerline{\hbox{\includegraphics[width=0.8in]{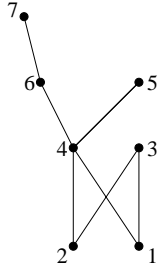}}}
\caption{A series parallel interval order and its preorder
extension.} \label{filt}
\end{center}
\end{figure}

For instance, the permutation associated with the poset depicted
in figure \ref{filt} is $2146753$. Indeed the filters of such a
poset are $R(1)=R(2)=\{3,4,5,6,7\}$, $R(3)=R(5)=R(7)=\emptyset$,
$R(4)=\{5,6,7\}$, $R(6)=\{7\}$ and then they are listed as
follows:
$$R(2)=R(1)\supset R(4)\supset R(6)\supset R(7)=R(5)=R(3).$$

\emph{Remark.} Given a permutation $\pi =a_1 \cdots a_n$ of length
$n$ and the partial order relation $R\in \mathcal{O}(n)$
associated with it, it is not difficult to observe that $R(a_j )$
is the set of the elements of $\pi$ greater than $a_j$ and
following $a_j$ in $\pi$. Analogously, $R^{-1}(a_j )$ is the set
of the elements of $\pi$ lesser than $a_j$ and preceding $a_j$ in
$\pi$ (see again figure \ref{filt}). In what follows, we will use
the notations $f_{\pi}(a_j )$ and $i_{\pi}(a_j )$ in place of
$R(a_j )$ and $R^{-1}(a_j )$ (respectively) when dealing with
permutations rather than partial order relations.

\subsection{The Tamari lattice and the weak Bruhat order on Av(312)}

According to \cite{K3}, it is possible to characterize the weak
Bruhat order on permutations using inversions. Recall that an
inversion of $\pi=a_1a_2...a_n$ is a pair $(a_i,a_j)$ such that
$i<j$ and $a_i>a_j$. Given two permutations of the same length
$\pi_1$ and $\pi_2$, it is $\pi_1 \leq_{b} \pi_2$ if and only if
the set $E(\pi_1)$ of inversions of $\pi_1$ is a subset of the set
$E(\pi_2)$ of inversion of $\pi_2$. The following simple
proposition provides the key ingredient to prove that the Tamari
lattice is isomorphic to the weak Bruhat order on $Av(312)$.

\begin{prop}
Let $\pi_1$ and $\pi_2$ be two permutations of length $n$. Then $
E(\pi_1) \subseteq E(\pi_2)$ if and only if, for every $1 \leq k
\leq n$, $f_{\pi_1}(k) \supseteq f_{\pi_2}(k)$.
\end{prop}
\emph{Proof.}\quad Suppose that, for every $k$, $f_{\pi_1}(k)
\supseteq f_{\pi_2}(k)$. If $(i,j) \in E(\pi_1)$, with $i>j$, then
$i \notin f_{\pi_1}(j)$ and then $i \notin f_{\pi_2}(j)$. This
implies that $(i,j) \in E(\pi_2)$.

Viceversa, suppose that $ E(\pi_1) \subseteq E(\pi_2)$. If $i \in
f_{\pi_2}(j)$, then $(i,j) \notin E(\pi_2)$, whence $(i,j) \notin
E(\pi_1)$. This implies that $i \in f_{\pi_1}(j)$.\cvd

\begin{cor}
The Tamari order is isomorphic to the weak Bruhat order restricted
to $Av_n(312)$.
\end{cor}

\subsection{The Dyck lattice and the strong Bruhat order on Av(312)}

For a given permutation $\pi$ of length $n$, define the vector
$max_{\pi}$ as follows: $max_{\pi}(k) = max\{\pi(i):i \leq k\}$.
According to \cite{BBFP}, we recall that, given two 312-avoiding
permutations $\pi_1$ and $\pi_2$ of length $n$, $\pi_1 \leq_{B}
\pi_2$ if and only if, for all $1 \leq k \leq n$, $max_{\pi_1}(k)
\leq max_{\pi_2}(k)$. For instance, considering the two
permutations $\pi_1=468753921, \pi_2=768543921\in Av_9(312)$, we
have $max_{\pi_1}=(4,6,8,8,8,8,9,9,9)$ and
$max_{\pi_2}=(7,7,8,8,8,8,9,9,9)$, whence $\pi_1 \leq_{B} \pi_2$.
Indeed, starting from $\pi_2$, we obtain $\pi_1$ by the following
reductions: $\pi_2=768543921 \rightarrow 768453921 \rightarrow
468753921=\pi_1$. Observe that, in the above sequence of
reductions, the permutation 768453921 is not 312-avoiding.

\bigskip

Given a permutation $\pi$ of length $n$, consider the set of its
\emph{consecutive noninversions}, i.e. the set of all $m \in
\{1,...,n\}$ such that either $m = n$ or the pair $(m,m+1)$ is a
noninversion of $\pi$ (that is $m$ appears before $m+1$ in $\pi$).
The following lemma provides a characterization of consecutive
noninversions in permutations avoiding $312$.

\begin{lemma}
Let $\pi \in Av_n(312)$, then the following properties hold:
\begin{enumerate}

\item[i)] $m$ is a consecutive noninversion of $\pi$ if and only
if either $m=n$ or $|i_{\pi}(m)| \neq |i_{\pi}(m+1)|$;

\item[ii)] if $j<k$, then $|i_{\pi}(j)| \leq |i_{\pi}(k)|$;

\item[iii)] if $k = |i_{\pi}(m)|$ and $m$ is a consecutive
noninversion of $\pi$, then $\pi(k+1) = m$;

\item[iv)] the consecutive noninversions of $\pi$ are those
elements of $\pi$ preceded only by lesser entries;

\item[v)] the set of all consecutive noninversions of $\pi$
coincides with the set of components of $max_{\pi}$. Moreover, the
index of the consecutive noninversion $m$ in $\pi$ coincides with
the index of the first occurrence of $m$ in $max_{\pi}$.
\end{enumerate}
\end{lemma}

\emph{Proof.}\quad \emph{i)} If $m\neq n$ is a consecutive
noninversion, then obviously $|i_{\pi}(m+1)|\geq |i_{\pi}(m)|-1$,
whence $|i_{\pi}(m)|\neq |i_{\pi}(m+1)|$. Viceversa, suppose that
$m\neq n$ is such that $|i_{\pi}(m)| \neq |i_{\pi}(m+1)|$; if
$(m+1,m)$ were an inversion of $\pi$, then there should be an
entry $j<m$ between $m+1$ and $m$, which is impossible since $\pi
\in Av_n (312)$.

\emph{ii)} This is an immediate consequence of the fact that $\pi
\in Av_n (312)$.

\emph {iii)} This is obvious when $m = n$. Otherwise, suppose that
$m$ is an element of $\pi$ having more than $i$ elements on its
left; then there would be at least one element $j$ such that $j$
precedes $m$ and $j>m$, and the three elements $j,m,m+1$ would
show a 312-pattern in $\pi$, a contradiction.

\emph{iv)} Observe that, if $m$ is a consecutive noninversion of
$\pi$, then, by \emph{iii)}, all elements of $\pi$ preceding $m$
are less than $m$. Viceversa, if $m$ is an entry of $\pi$ preceded
only by lesser elements, one cannot have $|i_{\pi}(m)| =
|i_{\pi}(m+1)|$, since in this case $m+1$ would precede $m$. Then
$m = n$ or $|i_{\pi}(m)| < |i_{\pi}(m+1)|$.

\emph {v)} This is a direct consequence of {\emph {iv}}).\cvd

\bigskip

\begin{prop}
Given $\pi_1 ,\pi_2 \in Av_n (312)$, the following conditions are
equivalent:
\begin{itemize}

\item[a)] for every $1\leq k\leq n$, $|i_{\pi_{1}}(k)|\geq
|i_{\pi_{2}}(k)|$;

\item[b)] for every $1\leq k\leq n$, $max_{\pi_{1}}(k)\leq
max_{\pi_{2}}(k)$.

\end{itemize}
\end{prop}
\emph{Proof.}\quad $a) \Rightarrow b)$ Let $j'=\pi_1 (j)$ be a
consecutive noninversion of $\pi_1$. Thanks to the above lemma,
item \emph{iii)}, all the elements of $\pi_1$ preceding $j'$ are
less than $j'$. We claim that there exists $k=\pi_2 (h)$ such that
$k>j'$, $h\leq j$ and all the elements of $\pi_2$ before $k$ are
less than $k$. Indeed, we have $j-1 = |i_{\pi_1}(j')| \geq
|i_{\pi_2}(j')| = |i_{\pi_2}(k)| = t $, where $k$ is the first
consecutive noninversion of $\pi_2$ on the left of $j'$ with
$k\geq j'$ (such a $k$ does indeed exist, as the reader can
immediately check). Then $k=\pi_2 (t+1)$ is the desired element of
$\pi_2$. As a consequence, we have that, if $j'$ is a consecutive
noninversion of $\pi_1$, then $max_{\pi_1}(j')\leq
max_{\pi_2}(j')$. From this we can immediately deduce the same
inequality for any $j'\leq n$.

$b) \Rightarrow a)$ Set $t = |i_{\pi_1}(j)|$, consider the
consecutive noninversion $j'$ of $\pi_1$ such that
$|i_{\pi_1}(j')| = t$ with $j' \geq j$ (this is the first
consecutive noninversion of $\pi_1$ on the left of $j$). From the
previous lemma we have $j' = max_{\pi_1}(t+1) \leq
max_{\pi_2}(t+1) = max_{\pi_2}(h)$, where $h \leq t+1$ is the
index of the first component of $max_{\pi_2}$, which is equal to
$max_{\pi_2}(t+1)$. Again from the previous lemma, we have that
$max_{\pi_2}(h) = k$ is a consecutive noninversion of $\pi_2$ and
$|i_{\pi_2}(k)| = h-1$. Finally, using in particular item
\emph{ii)} of the above lemma in the first two inequalities, we
have :
$$|i_{\pi_2}(j)| \leq |i_{\pi_2}(j')| \leq |i_{\pi_2}(k)| = h-1 \leq t = |i_{\pi_1}(j)|.$$

\cvd

\begin{cor}
The Dyck order is isomorphic to the strong Bruhat order restricted
to $Av_n(312)$.
\end{cor}

\section{Further works}

In the present work we have considered two well known Catalan
posets and we have proposed a unifying language to describe them
based on the notion of series parallel interval order. There are
of course several other poset structures which can be considered
on the objects of the Catalan family. Maybe the most famous one is
the \emph{Kreweras order} \cite{K}, which is naturally defined on
noncrossing partitions of a set of given cardinality by refining
the classical partial order on set partition. Other less classical
posets have been defined by Baril and Pallo on Dyck words (the
\emph{phagocyte lattice} \cite{BP1}) and on binary trees (the
\emph{pruning-grafting lattice} \cite{BP2}). It seems natural to
ask if series parallel interval orders can be used also to
describe these (and maybe other) Catalan posets. Unfortunately, we
have not been able to find an answer to such a question yet.

\end{document}